\documentclass[11pt,psamsfonts]{amsart}
\usepackage{amsmath}
\usepackage{amsthm}
\usepackage{amssymb}
\usepackage{amscd}
\usepackage{amsfonts}
\usepackage{amsbsy}
\usepackage{epsfig,afterpage}
\usepackage[dvips]{psfrag}
\newcommand{\R}{\ensuremath{\mathbb{R}}}

\newcommand{\Z}{\ensuremath{\mathbb{Z}}}

\newcommand{\e}{\varepsilon}

\newcommand{\dem}{\noindent {\it Proof: }}
\newcommand{\demP}{\noindent {\it Proof of Proposition \ref{prop1c}: }}
\newcommand{\cqd}{\hfill $\Box$}

\newtheorem{observation}{Remark}[section]

\newtheorem{definition}[observation]{Definition}

\newtheorem{proposition}[observation]{Proposition}
\newtheorem{theorem}[observation]{Theorem}
\newtheorem{lemma}[observation]{Lemma}

\def\Fix{\mathrm{Fix}}
\def\Im{\mathrm{Im}}
\def\Ker{\mathrm{Ker}}
\def\X{X_{H}}

\def\x{x}
\def\y{y}

\begin{document}

\title[Periodic orbits in reversible
hamiltonian systems] {Branching of periodic orbits in reversible
hamiltonian systems}

\author[ C.A. Buzzi, L.A. Roberto, M.A. Teixeira]
{Claudio A. Buzzi$^1$, Luci Any Roberto$^2$ and Marco A.
Teixeira$^3$}

\address{$^1$ Departamento de Matem\'atica, IBILCE--UNESP,
Rua C. Colombo, 2265, CEP 15054--000 S. J. Rio Preto, S\~ao Paulo,
Brazil.}

\address{$^2$ Departamento de Matem\'atica, IMECC--UNICAMP, CEP 13081--970, Campinas, S\~ao Paulo, Brazil.}

\address{$^3$ Departamento de Matem\'atica, IMECC--UNICAMP, CEP 13081--970, Campinas, S\~ao Paulo, Brazil.}

\email{buzzi@ibilce.unesp.br}

\email{lroberto@ime.unicamp.br}

\email{teixeira@ime.unicamp.br}

\date{}
\dedicatory{} \maketitle

\begin{abstract}
This paper deals with the dynamics of time-reversible Hamiltonian
vector fields with $2$ and $3$ degrees of freedom around an elliptic
equilibrium point in presence of symplectic involutions. The main
results discuss the existence of one-parameter families of
reversible periodic solutions terminating at the equilibrium. The
main techniques used are Birkhoff and Belitskii normal forms
combined with the Liapunov-Schmidt reduction.
\end{abstract}

%\part{Use this type of header for very long papers only}
% use lowercase except for proper names

\section{Introduction} % use lowercase except for proper names
\label{intro}

The resemblance of dynamics between reversible and Hamiltonian
contexts, probably first noticed by Poincar\'{e} and Birkhoff, has
caught much attention since the sixties of the twentieth century.
Since then many important results, e.g. KAM theory, Liapunov center
theorems, etc, holding in the Hamiltonian context have been carried
over to the reversible one (see \cite{KV,VV} and reference therein).

\smallskip

%The classical theorem due to Liapunov (see for example in Abraham
%and Marsden \cite{AM} p 498) establishes the existence of families
%of periodic orbits arising from an equilibrium point for real
%Hamiltonian systems under the following conditions: (i) the $2n$
%eigenvalues $\pm \lambda_1,\ \pm \lambda_2,\dots,\pm \lambda_n$ are
%distinct, (ii) $\lambda_1$ is purely imaginary and (iii)
%$\frac{\lambda_p}{\lambda_1}$ is non integer for $p=2,3,\dots,n$
%(non-resonance). Each such family of periodic orbits is
%parameterized by some real parameter $\sigma$ such that as
%$\sigma\rightarrow 0$, the periodic orbits tend towards the
%equilibrium while their periods tend towards
%$\frac{2\pi}{|\lambda_1|}$.

The concept of reversibility is linked with an involution $R$, i.
e., a map $R:\R^N\rightarrow\R^N$ such that $R\circ R=Id$. Let $X$
be a smooth vector field on $R^{N}$. The vector field is called
$R$--reversible if the following relation is satisfied
\[X(R(x))=-DR_x.X(x).\]
Reversibility means that $x(t)$ is a solution of $X$ if and only if
$Rx(-t)$ is also a solution. The set $Fix(R)=\{x\in\R^N:R(x)=x\}$
plays an important role in the reversible systems. We say that a
singular point $p$ is symmetric if $p\in Fix(R)$, and analogously we
say that an orbit $\gamma$ is symmetric if $R(\gamma)=\gamma$.

\smallskip
%One of characteristics properties of reversible systems is that
%generically  periodic orbits or invariant tori or minimal sets of
%such systems typically appear in one-parameter families.

Many dynamical systems that arise in the context of applications
possess robust structural properties, such as for instance
symmetries or Hamiltonian structure. In order to understand the
typical dynamics of such systems, their structure need to be taken
into account, leading one to study phenomena that are generic among
dynamical systems with the same structure. In the last decade there
has been a surging interest in the study of systems with
time-reversal symmetries (see \cite{Te} and \cite{Ha}). Symmetry
properties arise naturally and frequently in dynamical systems. In
recent years, a lot of attention has been devoted to understand and
use the interplay between dynamics and symmetry properties. It is
worthwhile to mention that one of the characteristic properties of
Hamiltonian and reversible systems is that minimal sets appear in
one-parameter families. So a number of natural questions can be
formulated, such as: (i) how do branches of such minimal sets
terminate or originate?; (ii) can one branch of minimal sets
bifurcate from another such branch?; (iii) how persistent is such
branching process when the original system is slightly perturbed?
Recently, there has been increased interest in the study of systems
with time-reversal symmetries and we refer \cite{LR} for a survey in
reversible systems and related problems.

\smallskip

Our main concern, in this article, is to find conditions for the
existence of one-parameter families of periodic orbits terminating
at the equilibrium.

\smallskip

We present  some relevant historical facts. In 1895 Liapunov
published his celebrated center theorem, see Abraham and Marsden
\cite{AM} p 498; This theorem, for analytic Hamiltonians with $n$
degrees of freedom, states that if the eigenfrequencies of the
linearized Hamiltonian are independent over $\Z$, near a stable
equilibrium point, then there exists $n$ families of periodic
solutions filling up smooth 2-dimensional manifolds going through
the equilibrium point. Devaney \cite{De} proved a time-reversible
version of the Liapunov center theorem. Recently this center theorem
has been generalized to equivariant systems, by Golubitsky, Krupa
and Lim \cite{GKL} in the time-reversible case, and by Montaldi,
Roberts and Stewart \cite{MRS} in the Hamiltonian case. We recall
that in \cite{GKL} the Devaney's theorem was extended and some extra
symmetries were considered. Contrasting Devaney's geometrical
approach, they used Liapunov-Schmidt reduction, adapting an
alternative proof of the reversible Liapunov center theorem given by
Vanderbauwhede \cite{Va}. In \cite{MRS} the existence of families of
periodic orbits around an elliptic semi-simple equilibrium is
analyzed. Systems with symmetry, including time-reversal symmetry,
which are anti-symplectic are studied. Their approach is a
continuation of the work of Vanderbauwhede, in \cite{Va}, where the
families of periodic solutions correspond bijectively to solutions
of a variational problem.

\smallskip

%The reversible version of the Liapunov Theorem due to Devaney (see
%\cite{De}) says that if $X$ is a $C^2$ $R$--reversible vector field
%in a neighborhood of a symmetric singular point $p$, its linear part
%posses eigenvalues
%$\pm\lambda=\pm\lambda_1,\pm\lambda_2,\dots,\pm\lambda_n$, $\lambda$
%is purely imaginary, and no other eigenvalue is equal to $i\lambda
%m$ for any integer $m$, then there exists a $C^2$, two--dimensional,
%invariant manifold $M^\lambda$ containing $p$ with the property that
%$M^\lambda$ consists of a nested, one--parameter family of symmetric
%periodic orbits. Moreover, the periods of the closed orbits tend to
%$2\pi/|\lambda|$ as the initial conditions tend to $p$.

Recently Buzzi and Teixeira in \cite{BT} have analyzed  the dynamics
of time-reversible Hamiltonian vector fields with 2 degrees of
freedom around an elliptic equilibrium point in presence of $1 :-1$
resonance. Such systems appear generically inside a class of
Hamiltonian vector fields in which the symplectic structure is
assumed to have some symmetric properties. Roughly speaking, the
main result says that under certain conditions the original
Hamiltonian $H$ is formally equivalent to another Hamiltonian
$\widetilde{H}$ such that the corresponding Hamiltonian vector field
$X_{\widetilde{H}}$ has two Liapunov families of symmetric periodic
solutions terminating at the equilibrium. It is worth while to say
that all the systems considered there have been derived from the
expression of Birkhoff normal form.

\smallskip

In this paper we address the problem to systems with 2 and 3 degrees
of freedom. Physical models of such systems were exhibited in
\cite{CK,He}. As usual the main proofs are based on a combined use
of normal form theory and the Liapunov-Schmidt Reduction. It is
important to mention that our results concerning the existence of
Liapunov families generalize those in \cite{BT}. As a matter of fact
we deal with $C^\infty$ or $C^{\omega}$.
\smallskip

We begin in Section 2 with an introduction of the terminology and
basic concepts for the formulation of our results. In Section 3 the
Belitskii normal form is discussed. In Section 4 the
Liapunov-Schmidt reduction is presented. In Section 5 the usefulness
of Birkhoff normal form in our approach is pointed out. In Section 6
we study the Hamiltonian with 2 degrees of freedom denoted by
$\Omega^0$, and we denote by $\Omega^0_B$ the set of vector fields
in $\Omega^0$ that satisfy the Birkhoff Condition and by
$\Omega^0_\omega$ the vector fields in $\Omega^0$ that are analytic.
We generalize some results presented in \cite{BT} by proving Theorem
A. That result says that there exists an open set
$\mathcal{U}^0\subset\Omega^0_B$ (respec. $\Omega^0_\omega$), in the
$C^\infty$--topology, such that (a) $\mathcal{U}^0$ is determined by
the $3$--jet of the vector fields; and (b) each $X\in\mathcal{U}^0$
possesses two $1$--parameter families of periodic solutions
terminating at the equilibrium. In Section 7 we study the
Hamiltonian with 3 degrees of freedom, and we prove Theorems B and
C. In Theorem B we consider the involution associated to the system
satisfying $dim (Fix(R))=2$, and in Theorem C satisfying $dim
(Fix(R))=4$. We denote these spaces of reversible Hamiltonian vector
fields by $\Omega^1$ and $\Omega^2$, respectively. Again
$\Omega^2_B$ is the set of vector fields in $\Omega^2$ that satisfy
the Birkhoff Condition and $\Omega^2_\omega$ is the set of vector
fields in $\Omega^2$ that are analytic. The conclusions are the
following: In Theorem B there exists an open set
$\mathcal{U}^1\subset\Omega^1$, in the $C^\infty$--topology, such
that (a) $\mathcal{U}^1$ is determined by the $2$--jet of the vector
fields, and (b) for each $X\in\mathcal{U}^1$ there is no periodic
orbit arbitrarily close to the equilibrium. In Theorem C there
exists an open set $\mathcal{U}^2\subset\Omega^2_B$ (respec.
$\Omega^2_\omega$), in the $C^\infty$--topology, such that (a)
$\mathcal{U}^2$ is determined by the $3$--jet of the vector fields,
and (b) each $X\in\mathcal{U}^2$ has infinitely many one--parameter
family of periodic solutions terminating at an equilibrium with the
periods tending to $2\pi/\alpha$. In Section 8 we present an example
that satisfies the hypotheses of the Theorem A and commented on that
is possible to accomplish the vector fields of Theorem C.

\section{Preliminaries}

\smallskip

Now we introduce some of the terminology and basic concepts for the
formulation of our results.

\smallskip

We consider (germs of) smooth functions
$H:\mathbb{R}^{2n},0\rightarrow\mathbb{R}$ having the origin as an
equilibrium point. The corresponding Hamiltonian vector field, to be
denoted by $\X$, has the origin as an equilibrium or singular point.
We recall that $dH=\omega(\X,\cdot)$, where $\omega=dx_1\wedge
dy_1+dx_2\wedge dy_2+\cdots+dx_n\wedge dy_n$ denotes the standard
2-form on $\mathbb{R}^{2n}$. In coordinates $\X$ is expressed as:
$$\dot{x_i}=\frac{\partial H}{\partial
y_i},\,\,\,\dot{y_i}=-\frac{\partial H}{\partial
x_i};\,\,\,i=1,\cdots,n.$$ In $\mathbb{R}^6$ we have
$$\left(\begin{array}{c}\dot{x_1}\\\dot{y_1}\\ \vdots \\ \dot{x_3}\\\dot{y_3}
\end{array}\right)=
\left(\begin{array}{cccccc}
0 &1& 0& 0& 0& 0\\
-1& 0& 0& 0& 0& 0\\
0& 0& 0& 1& 0& 0\\
0& 0& -1& 0& 0& 0\\
0& 0& 0& 0& 0& 1\\
0& 0& 0& 0& -1&0\end{array}\right) \left(\begin{array}{c}
\frac{\partial H}{\partial x_1}\\
\frac{\partial H}{\partial y_1}
\\\ \vdots \\\frac{\partial H}{\partial x_3}\\
\frac{\partial H}{\partial y_3}\end{array}\right).$$

\smallskip

Here,
$$J=\left(\begin{array}{cccccc}0 &1& 0& 0& 0& 0\\-1& 0& 0& 0& 0& 0\\0& 0& 0& 1& 0& 0\\0& 0& -1& 0& 0&
0\\0& 0& 0& 0& 0& 1\\0& 0& 0& 0& -1& 0\end{array}\right)$$ is the
symplectic structure  associated with the 2-form $\omega$ given
above.

\smallskip

We say that an involution is \textit{symplectic} when it satisfies
the equation $\omega(DR_p(v_p),DR_p(w_p)) = \omega(v_p,w_p)$. If the
involution $R$ is linear then this definition is equivalent to
$JR=R^T J$, where $J$ is the symplectic structure and $R^T$ is the
transpose matrix of $R$.

\smallskip

The next proposition exhibits normal forms for linear symplectic
involutions on $\R^6$.

\smallskip

\begin{proposition}\label{prop1c}Given the symplectic structure $\omega$ and an
involution $R$ there exists a symplectic change of coordinates that
transforms $R$ in one of the following normal forms
\begin{enumerate}
\item $R_0=Id,$
\item $R_0(x_1,y_1,x_2,y_2,x_3,y_3)=(x_1,y_1,x_2,y_2,-x_3,-y_3),$
\item $R_0(x_1,y_1,x_2,y_2,x_3,y_3)=(x_1,y_1,-x_2,-y_2,-x_3,-y_3),$
\item $R_0=-Id.$
\end{enumerate}
\end{proposition}

\smallskip

Before giving the proof we observe that the mapping
$\psi=(1/2)(R+L)$, where $L=DR(0)$, is a symplectic conjugacy
between $R$ and $L$, i. e., $R\circ\psi=\psi\circ L.$ So we may and
do assume, without loss of generality, that the involution $R$ is
linear.

\smallskip

\begin{lemma}\label{lema2c}
If $R$ is a linear symplectic involution, then  we have that
$\mathbb{R}^6=\Fix(R)\oplus\Fix(-R)$ and
$\omega(\Fix(R),\Fix(-R))=0.$
\end{lemma}
\dem For every $u\in\mathbb{R}^6$, we can write
$u=((u+R(u))/2)+((u-R(u))/2)$. Notice that $(u+R(u))/2\in\Fix(R)$
and $(u-R(u))/2\in\Fix(-R).$ Now, let $u\in\Fix(R)$ and
$v\in\Fix(-R)$, so we have that $\omega(u,v)=\omega(R(u),-R(v))$. By
using that $R$ is symplectic and $R$ is linear, we have that
$\omega(R(u),R(v))=\omega(u,v)$. So $-\omega(u,v)=\omega(u,v)$, and
we have proved that $\omega(\Fix(R),\Fix(-R))$ $=0.$ \cqd

\smallskip

A linear subspace $U\in\mathbb{R}^6$ is {\it symplectic} if $\omega$
is {\it non-degenerate} in $U$, i. e, if $\omega(u,v)=0$ for all
$u\in U$ then $v=0$.

\smallskip

\begin{lemma}\label{lema3c} $\Fix(R)$ and $\Fix(-R)$ are symplectic
subspaces.\end{lemma}

\smallskip

\dem Suppose  $u\in\Fix(R)$ and $u\neq 0$ such that
$\omega(u,\Fix(R))=0$. By using Lemma \ref{lema2c}, we have
$\omega(\Fix(R),\Fix(-R))=0$, so $\omega(u,\Fix(-R))=0$. Again by
Lemma \ref{lema2c} ($\mathbb{R}^6=\Fix(R)\oplus\Fix(-R)$) we have
$\omega(u,\mathbb{R}^6)=0$ and so $\omega$ is degenerate in
$\mathbb{R}^6$ which is not true. Then $\Fix(R)$ is a symplectic
subspace. The proof for $\Fix(-R)$ is analogous. \cqd

\smallskip

\demP  Let $R:\mathbb{R}^6\rightarrow\mathbb{R}^6$ be a linear
involution and $\omega$ be a fixed symplectic structure. From Lemma
\ref{lema2c}, $\mathbb{R}^6=\Fix(R)\oplus\Fix(-R)$ and as $\Fix(R)$
is a symplectic subspace, then $\dim\Fix(R)=0,2,4,$ or $6$.

\smallskip

\begin{itemize}
\item if $\dim\Fix(R)=0$, then we can find a coordinate system, using Darboux Theorem \cite{HM}, such that $R_0=-Id$;

\smallskip

\item if $\dim\Fix(R)=6$, then we can find a coordinate system, using Darboux Theorem \cite{HM}, such that $R_0=Id$;

\smallskip

\item if $\dim\Fix(R)=4$, we consider the bases
$\beta_1=\{e_1,e_2,e_3,e_4\}$ for $\Fix (R)$ and
$\beta_2=\{f_1,f_2\}$ for $\Fix(-R)$. So
$\beta=\{e_1,e_2,e_3,e_4,f_1,f_2\}$ is a basis for $\mathbb{R}^6$.
Let us show that $\beta$ can be chosen such that
$[\omega]_{\beta}=J$ and
$[R]_{\beta}=R_0=\left(\begin{array}{cccccc}1&0&0&0&0&0\\0&1&0&0&0&0\\0&0&1&0&0&0\\0&0&0&1&0&0\\0&0&0&0&-1&0\\0&0&0&0&0&-1\end{array}\right)$.
Here $[\omega]_{\beta}$ means the matrix of $\omega$ with respect to
the basis $\beta$.

\smallskip

Note that $\omega(e_i,e_i)=0$ and $\omega(f_j,f_j)=0$, $i=1,2,3,4$
and $j=1,2$. By Lemma \ref{lema2c} $\omega(e_i,f_j)=0$, $i=1,2,3,4$
and $j=1,2$. And as $\omega$ is alternating, then
$\omega(f_1,f_2)=1$ and $\omega(f_2,f_1)=-1$.

\smallskip

Define $\omega(e_i,e_j)$ for $i\neq j$. From Darboux's Theorem there
exists a coordinate system around $0$ such that $\omega|_{\beta_1}$
in this coordinate system is the symplectic structure $J$. \item if
$\dim \Fix(R)=2$, in the same way as above, we get
$$\\R_0=[R]_{\beta}=\left(\begin{array}{cccccc}1&0&0&0&0&0\\0&1&0&0&0&0\\0&0&-1&0&0&0\\0&0&0&-1&0&0\\0&0&0&0&-1&0\\0&0&0&0&0&-1\end{array}\right).$$
\cqd
\end{itemize}

\medskip

Using the previous proposition we consider the following cases:

\smallskip

\begin{itemize}\item[]$6:2$--\textbf{Case}:
$R_1(x_1,y_1,x_2,y_2,x_3,y_3)=(x_1,$ $y_1,$ $-x_2,$ $-y_2,$ $-x_3,$
$-y_3)$,

\medskip

\item[]$6:4$--\textbf{Case}:
$R_2(x_1,y_1,x_2,y_2,x_3,y_3)=(x_1,y_1,x_2,y_2,$ $-x_3,$ $-y_3).$

\end{itemize}

\medskip

\subsection{Linear part of a $R_j$--reversible Hamiltonian vector field in $\mathbb{R}^6$}

Denote by $\Omega^j$ the space of all $R_j$--reversible Hamiltonian
vector field, $X_{H_j}$, in $\mathbb{R}^6$ with 3-degrees freedom
where $H_j$ is the associate Hamiltonian and $j=1,2$. Fix the
coordinate system $(x_1,y_{1},x_{2},y_{2},x_{3},y_{3})\in
(\mathbb{R}^6,0)$. We endow $\Omega^j$ with the
$C^{\infty}$--topology.

\smallskip

The symplectic structure given by $J$ is:
$$J=\left(\begin{array}{cccccc}0 &1& 0& 0& 0& 0\\-1& 0& 0& 0& 0& 0\\0& 0& 0& 1& 0& 0\\0& 0& -1& 0& 0&
0\\0& 0& 0& 0& 0& 1\\0& 0& 0& 0& -1& 0\end{array}\right).$$ Observe
that the involution $R_j$ is symplectic, i.e, $J.R_j - R_j^T. J=0$,
$j=1,2.$

\smallskip

As the involution is symplectic, then the vector field is
$R_j-$reversible if and only if the Hamiltonian function $H_j$ is
$R_j-$anti-invariant, $j=1,2$. This is equivalent to say that
$H_j\circ R_j=-H_j.$ (See \cite{BT})

\smallskip

Define the polynomial function with constant coefficients
$a_k\in\R$:

\smallskip

\noindent $H_j(x_1,y_{1},x_{2},y_{2},x_{3},y_{3}) =
    a_{01}x_1^2 + a_{02}x_1y_{1} + a_{03}x_1x_{2} + a_{04}x_1y_{2} + a_{05}x_1x_{3} + a_{06}x_1y_{3} +
      a_{07}y_{1}^2 + a_{08}y_{1}x_{2} + a_{09}y_{1}y_{2} + a_{10}y_{1}x_{3} + a_{11}y_{1}y_{3} + a_{12}x_{2}^2 +
      a_{13}x_{2}y_{2} + a_{14}x_{2}x_{3} + a_{15}x_{2}y_{3} + a_{16}y_{2}^2 + a_{17}y_{2}x_{3} + a_{18}y_{2}y_{3} +
      a_{19}x_{3}^2 + a_{20}x_{3}y_{3} + a_{21}y_{3}^2 + h.o.t..$

\smallskip

First of all we impose the $R_j$--reversibility on our Hamiltonian
system, $j=1,2$. For each case we have:

\smallskip

\begin{itemize}
\item[a)]\textbf{Case} $6:2$

\smallskip
From the reversibility condition, $H_1\circ R_1=-H_1$, and
$$R_1=\left( \begin{array}{cccccc} 1&0&0&0&0&0 \\
0&1&0&0&0&0\\0&0&-1&0&0&0\\0&0&0&-1&0&0\\0&0&0&0&-1&0\\0&0&0&0&0&-1
\end{array}\right),$$ we obtain
$$\begin{array}{ll}
H_1= & a_{03}x_1x_{2} + a_{04}x_1y_{2} + a_{05}x_1x_{3} + a_{06}x_1y_{3}+\\
&a_{08}x_{2}y_{1}+a_{09}y_{1}y_{2}+a_{10}x_{3}y_{1}+a_{11}y_{1}y_{3}
+ h.o.t..
\end{array}$$
Then, the linear part of Hamiltonian vector field $X_{H_1}$ is
$$A_1=\left(\begin{array}{cccccc}0& 0& a& b& c& d\\0& 0& e& f& g& h\\-f& b& 0& 0& 0& 0\\e& -a& 0& 0& 0& \
0\\-h& d& 0& 0& 0& 0\\g& -c& 0& 0& 0& 0
\end{array}\right).$$
Just to simplify the notation we replace $a_{03},$ $a_{04},$
$a_{05},$ $a_{06},$ $a_{08},$ $a_{09},$ $a_{10}$, $a_{11}$ by $a,$
$b,$ $c,$ $d,$ $-e$, $-f,$ $-g,$ $-h$, respectively. Note that $A_1$
is $R_1-$reversible (i. e, $R_1.A_1 + A_1.R_1=0$). The eigenvalues
of $A_1$ are $\{0,0, \pm\sqrt{be-af+dg-ch},
\pm\sqrt{be-af+dg-ch}\}$. We restrict our attention to those systems
satisfying  the inequality:
\begin{equation} \label{naosei}
be-af+dg-ch <0.
\end{equation}
The case when $be-af+dg-ch >0$ will not be considered because the
center manifold of the equilibrium has dimension two with double
zero eigenvalue. We shall use the Jordan canonical form from $A_1$.
So we stay, for while, away from the original symplectic structure.
We call $\alpha=\sqrt{-be+af-dg+ch}$, and so the transformation
matrix is

$$P_1=\left(\begin{array}{cccccc}
0&0&\frac{-d}{dg-ch}\alpha&0&\frac{-c}{dg-ch}\alpha&0\\
0&0&\frac{-h}{dg-ch}\alpha&0&\frac{-g}{dg-ch}\alpha&0\\
\frac{df-bh}{be-af}&\frac{cf-bg}{be-af}&0&\frac{-df+bh}{dg-ch}&0&\frac{-cf+bg}{dg-ch}\\
\frac{-de+ah}{be-af}&\frac{-ce+ag}{be-af}&0&\frac{de-ah}{dg-ch}&0&\frac{ce-ag}{dg-ch}\\
0& 1& 0& 0& & 1\\1& 0& 0& 1& 0& 0
\end{array}\right).$$ So
$$\widehat{A_1}=P_1^{-1}.A_1.P_1=
\left(\begin{array}{cccccc}
0&0&0&0&0&0\\
0&0&0&0&0&0\\
0&0&0&\alpha&0&0\\
0&0&-\alpha&0&0&0\\
0&0&0&0&0&\alpha\\
0&0&0&0&-\alpha&0
\end{array}\right),$$ where $P_1^{-1}$ is the inverse matrix of the matrix $P_1$. Moreover, in this
way, $\widehat{R_1}=P_1^{-1}.R_1.P$ takes the form
$$\widehat{R_1}=
\left(\begin{array}{cccccc}
-1& 0& 0& 0& 0& 0\\
0& -1& 0& 0& 0& 0\\
0& 0& 1& 0& 0& 0\\
0& 0& 0& -1& 0& 0\\
0& 0& 0& 0& 1& 0\\
0& 0& 0& 0& 0& -1
\end{array}\right).$$

\item[b)]\textbf{Case} $6:4$

\smallskip

We proceed in the same way as in the previous case. The involution
is
$$R_2=\left( \begin{array}{cccccc}
1&0&0&0&0&0 \\
0&1&0&0&0&0\\
0&0&1&0&0&0\\
0&0&0&1&0&0\\
0&0&0&0&-1&0\\
0&0&0&0&0&-1
\end{array}\right)$$
and the Hamiltonian function in this case takes the form:

$$\begin{array}{ll}H_2=&
a_{05}x_1x_{3}+a_{14}x_{2}x_{3}+a_{10}x_{3}y_{1}+a_{17}x_{3}y_{2}+\\
&a_{06}x_1y_{3}+a_{15}x_{2}y_{3}+a_{11}y_{1}y_{3}+a_{18}y_{2}y_{3} +
h.o.t..
\end{array}$$

\smallskip

Then, the linear part of Hamiltonian vector field $X_{H_2}$ is
expressed by:
$$A_2=\left(\begin{array}{cccccc}
0& 0& 0& 0& a& b\\
0& 0& 0& 0& c& d\\
0& 0& 0& 0& e& f\\
0& 0& 0& 0& g& h\\
-d& b& -h& f& 0& 0\\
c& -a& g& -e& 0&0
\end{array}\right).$$ Again we change the notation. The eigenvalues
of $A_2$ are given by $\{0,0, \pm\sqrt{bc-ad+fg-eh}, $
$\pm\sqrt{bc-ad+fg-eh}\}$. We consider the case
\begin{equation} \label{naosei2}
bc-ad+fg-eh <0.
\end{equation}

\smallskip

We call $\alpha=\sqrt{-bc+ad-fg+eh}$ and consider the transformation
matrix

\smallskip

$$P_2=\left(\begin{array}{cccccc}
\frac{be-af}{bc-ad}& \frac{-bg+ah}{bc-ad}&0&\frac{-b}{\alpha}&0&\frac{-a}{\alpha}\\
\frac{de-cf}{bc-ad}&\frac{-dg+ch}{bc-ad}&0&\frac{-d}{\alpha}&0&\frac{-c}{\alpha}\\
0&1&0&\frac{-f}{\alpha}&0&\frac{-e}{\alpha}\\
1&0&0&\frac{-h}{\alpha}&0&\frac{-g}{\alpha}\\
0& 0& 0& 0& 1& 0 \\
0& 0& 1& 0& 0& 0
\end{array}\right),$$
and the Jordan canonical form of $A_2$ is:
$$\widehat{A_2}=P_2^{-1}.A_2.P_2=
\left(\begin{array}{cccccc}
0&0&0&0&0&0\\
0&0&0&0&0&0\\
0&0&0&\alpha&0&0\\
0&0&-\alpha&0&0&0\\
0&0&0&0&0&\alpha\\
0&0&0&0&-\alpha&0
\end{array}\right).$$

\smallskip

Moreover, in this way, $\widehat{R_2}=P_2^{-1}.R_2.P_2$ takes the
form
$$\widehat{R_2}=\left(\begin{array}{cccccc}
1& 0& 0& 0& 0& 0\\
0& 1& 0& 0& 0& 0\\
0& 0& -1& 0& 0& 0\\
0& 0& 0& 1& 0& 0\\
0& 0& 0& 0& -1& 0\\
0& 0& 0& 0& 0& 1
\end{array}\right).$$
\end{itemize}

\medskip

\section{Belitskii normal form}

\smallskip

In this section we present the Belitskii Normal Form. When a vector
field is in this normal form we can write explicitly the resultant
equation of Liapunov--Schmidt reduction.

\smallskip

Consider a formal vector field expressed by
$$\hat{X}(\x)=A\x+\sum_{k\geq 2}X^{(k)}(\x)$$ where $X^{(k)}$
is the homogeneous part of degree $k$. Let us look for a ``simple"
form of the formal vector field $\hat{Y}=\hat{\phi}*\hat{X}$ by
means of formal transformation
$$\hat{\phi}=\x+\sum_{k}^{\infty}\phi^{(k)}(\x).$$

\smallskip

The proof of the next theorem is in \cite{KI}.

\smallskip

\begin{theorem}
Given a formal vector field $$\hat{X}(\x)=A\x+\sum_{k\geq
2}X^{(k)}(\x),$$ there is a formal transformation $\hat{\phi}(\x)=\x
+\dots$ bringing $\hat{X}$ to the form $(\hat{\phi}*X)(\x)=Ax +
h(\x)$ where $h$ is a formal vector field with zero linear part
commuting with $A^T$, i.e
$$A^T h(\x)=h'(\x)A^T \x,$$ where $A^T$ is the transposed matrix and
$h^\prime$ is the derivative of $h$.
\end{theorem}

\smallskip

Here we call the normal form $(\hat{\phi}_* X)(\x)=Ax + h(\x)$ the
Belitskii normal form. By abuse of the terminology, call $\X=A+h.$

\medskip

\section{Liapunov--Schmidt reduction}

\smallskip

In this section we recall the main features of the Liapunov--Schmidt
reduction. As a matter of fact, we adapt the setting presented in
\cite{CH,WA} to our approach. In this way consider the
$R$-reversible system expressed by

\smallskip

\begin{equation}\label{sistem}
\dot{\x}=X_{H}( \x);\,\, \x\in\mathbb{R}^6
\end{equation}
satisfying $X_{H}(R\x)=-RX_{H}(\x)$ with $R$ a linear involution in
$\mathbb{R}^6$. Assume that $\X(0)=0$ and consider
\begin{equation} A=D_1 \X(0),
\end{equation}
the Jacobian matrix of $\X$ in the origin.

\smallskip

In our case the linear part of vector field has the following
eigenvalues: $0$ with the algebraic and geometric multiplicity $2$,
and  $\pm \alpha i,$ also with algebraic and geometric multiplicity
$2$, $\alpha\in \,\mathbb{R}.$  Performing a time rescaling we may
take $\alpha =1$. We write the real form of the linear part of the
vector field $X_j$:
$$A=\left(
\begin{array}{cccccc}0&0&0&0&0&0\\0&0&0&0&0&0\\0&0&0&1&0&0\\0&0&-1&0&0&0
\\0&0&0&0&0&1\\0&0&0&0&-1&0 \end{array}\right).$$

\smallskip

Let $C^{0}_{2\pi}$ the Banach space of de $2\pi-$periodic continuous
mappings $\x:\mathbb{R}\rightarrow\mathbb{R}^6$ and $C^1_{2\pi}$ the
corresponding $C^1-$subspace. We define an inner product on
$C^{0}_{2\pi}$ by
$$(\x_1,\x_{2})=\frac{1}{2\pi}\int_{0}^{2\pi}<\x_1(t),\x_{2}(t)>dt$$ where
$<\cdot,\cdot>$ denotes an inner product in $\mathbb{R}^6$.

\smallskip

The main aim is to find all small periodic solutions of
(\ref{sistem}) with period near $2\pi$.

\smallskip

Define the map $F:C^{1}_{2\pi}\times \mathbb{R}\rightarrow
C^{0}_{2\pi}$  by
$$F(\x,\sigma)(t)=(1+\sigma)\dot{\x}(t)-\X(\x(t)).$$

\smallskip

Note that if $(\x_0,\sigma_0)\in C^{1}_{2\pi}\times\mathbb{R}$ is
such that
\begin{equation}\label{solution}F(\x_0,\sigma_0)=0,\end{equation}
then ${\tilde\x}(t):=\x_0((1+\sigma_0)t)$ is a
$2\pi/(1+\sigma_0)-$periodic solution of (\ref{sistem}).

\smallskip

Our task now is to find the zeroes of $F$. Clearly,
$(\x_0,\sigma_0)=(0,0)$ is one solution of $F(\x_0,\sigma_0)=0$. Let
$L:=D_{\x}F(0,0):C^{1}_{2\pi}\rightarrow C^{0}_{2\pi};$ explicitly
$L$ is given by
$$L\x(t)=\dot\x(t)-A\x(t).$$

\smallskip

Consider the unique (S-N)-decomposition of $A$, $A=S+N.$ Recall that
in our case $A$ is semi-simple, i. e, $A=S$. Define the subspace
$\mathcal{N}$ of $C^{1}_{2\pi}$ as
$$\begin{array}{ll}\mathcal{N}=&\{q; \dot{q}(t)=S q(t)\}=\\ &\{q;\,\, q(t)=exp(tS)x;\,\,x\in\mathbb{R}^6\}.\end{array}$$

\smallskip

Observe that $\mathcal{N}\subset C^{1}_{2\pi}$ and a basis for the
solutions of $\dot{q}=Sq$ is given by the set $\{(1,0,0,0,0,0),$
$(0,1,0,0,0,0),$ $(0,0,\cos(t),\sin(t),0,0),$
$(0,0,-\sin(t),\cos(t),$ $0,0),$ $(0,0,0,0,\cos(t),\sin(t)),$
$(0,0,0,0,-\sin(t),\cos(t))\}$.

\smallskip

In order to study certain properties of the operator $L$ we
introduce $\mathcal{N}\subset C^{1}_{2\pi}$ and the following
definitions and notations.

\smallskip

We will put the solution of $F(x_0,\sigma_0)=0$ in one-to-one
correspondence with the solutions of an appropriate equation in
$\mathcal{N}.$ Define the subspaces $$X_1=\{\x\in
C^{1}_{2\pi}:\,\,(\x,\mathcal{N})=0\}$$ and $$Y_1=\{\y\in
C^{0}_{2\pi}:(\y,\mathcal{N})=0\}$$ as the orthogonal complements of
$\mathcal{N}$ in $C^{1}_{2\pi}$ and $C^{0}_{2\pi}$, respectively.

\smallskip

Let $(q_1,q_2,q_3,q_4,q_5,q_6)$ with $q_i=exp(tS)u_i$ where $u_i$,
$i=1,...,6$, is a basis for $\mathbb{R}^6.$ Then we define a
projection $$\mathcal{P}:C^{0}_{2\pi}\rightarrow C^{0}_{2\pi}$$ by
$$\mathcal{P}=\sum_{i=1}^{6}q_i^{*}(\cdot)q_i\,\,\in\,\,\mathcal{L}(C^{0}_{2\pi})$$
with $q_i^{*}(\x)=(q_i,\x).$

\smallskip

We have $\Im(\mathcal{P})=\mathcal{N}$ and $\Ker(\mathcal{P})=Y_1.$
Hence,
$$C^{1}_{2\pi}=X_1\oplus\mathcal{N},\,\,\,C^{0}_{2\pi}=Y_1\oplus\mathcal{N}.$$

\smallskip

Now we consider
$$F(\x,\sigma)=F(q+\x_1,\sigma)=:\hat{F}(q,\x_1,\sigma);\,\,q\in\mathcal{N},\,\,\x_1\in X_1.$$

\smallskip

The proof of next result can be found in \cite{JH}.

\smallskip

\begin{lemma}({\bf Fredholm's Alternative}) Let $A(t)$ be a matrix
in $C^{0}_{T}$ and let $f$ be in $C_T$. Here $C^0_T$ is the space of
the matrices with entries continuous and $T$--periodic, and $C_T$ is
the set of $T$-periodic maps from $\R$ to $\R^n$. Then the equation
$\dot{\x}=A(t)\x+f(t)$ has a solution in $C_T$ if, and only if,
$$\int_{0}^{T}<\y(t),g(t)>dt=0$$ for all solution $\y$ of the
adjoint equation $$\dot{\y}=-\y A(t)$$ such that $\y^{t}\in C_T$.
\end{lemma}

\smallskip

As $L(\mathcal{N})\subset\mathcal{N}$ this lemma implies the
following:

\smallskip

\begin{lemma}\label{lema-bijecao}
The mapping $\hat{L}:=L|_{X_1}:X_1\rightarrow Y_1$ is bijective.
\end{lemma}

\smallskip

Let us study the solutions of $\hat{F}(q,x_1,\sigma)=0.$ These
solutions are equivalent to the solutions of the system
$$\begin{array}{r}
(I-\mathcal{P})\circ\hat{F}(q,x_1,\sigma)=0,\\
\mathcal{P}\circ\hat{F}(q,x_1,\sigma)=0.
\end{array}$$

\smallskip

With Lemma \ref{lema-bijecao} and the Implicit Function Theorem we
can solve the first equation as $\x_1=\x_1^{*}(q,\sigma).$ Then,
(\ref{solution}) is reduced to
$$\tilde{F}(q,\sigma):=\mathcal{P}\circ\hat{F}(q,\x_1^{*}(q,\sigma),\sigma)=0.$$

\smallskip

This equation is solved if, and only if,
$$q_i^{*}(\hat{F}(q,\x_1^{*}(q,\sigma),\sigma)=0,\,\,i=1,\cdots,6.$$

\smallskip

Notice that $(u,\sigma)$ is a solution of (\ref{solution}) provided
that
\begin{equation}\label{eqB}B(u,\sigma)=0\end{equation} with
$B:\mathcal{N}\times\mathbb{R}\rightarrow\mathbb{R}^6$ defined by

\smallskip

$$B(u,\sigma):=\frac{1}{2\pi}\int_{0}^{2\pi}exp(-tS)F(\x^{*}(u,\sigma),\sigma)dt$$
and
$$\x^{*}(u,\sigma):=exp(tS)u+\x^{*}_{1}(exp(tS)u,\sigma).$$

\smallskip

Let us present some properties of the mapping $B$.

\smallskip

The proof of next lemma can be found in \cite{KV}.

\smallskip

\begin{lemma}
The following relations hold:
\begin{itemize}
\item[i)] $s_{\phi}B(u,\sigma)=B(s_{\phi}u,\sigma);$
\item[ii)] $RB(u,\sigma)=-B(Ru,\sigma),$
where $s_{\phi}$ is the $S^1-$action in $\mathbb{R}^6$ defined by
$s_{\phi}u=exp(-\phi S_{0})u.$
\end{itemize}
\end{lemma}

\smallskip

Observe that under the condition $i)$ the mapping $B$ is $S^1-$
equivariant whereas condition $ii)$ states that the mapping $B$ is
$R-$anti-equivariant, i.e., $B$ inherits the anti-symmetric
properties of $\X.$

\smallskip

Assume that (\ref{sistem}) is in Belitskii normal form truncated at
the order $p$. So $\X(\x)=A\x+h(\x)+r(\x)$ where
$r(\x)=\mathcal{O}(\parallel\x \parallel^{p+1}).$ The proof of next
result is in \cite{WA}.

\smallskip

\begin{theorem}\label{teowagen}
The following relations hold:
\begin{enumerate}
\item[i)] $\x^*(u,\sigma)=exp(tS)u+\mathcal{O}(\parallel\x
\parallel^{p+1})$,
\item[ii)]
$B(u,\sigma)=(1+\sigma)Su-Au-h(u)+\mathcal{O}(\parallel\x
\parallel^{p+1})$
for $\sigma$ near the origin.
\end{enumerate}
\end{theorem}

\smallskip

If $(u,\sigma)$ is a solution of (\ref{eqB}) then $x=x^*(u,\sigma)$
corresponds to a  $2\pi/(1+\sigma)$-periodic solution of
(\ref{solution}).

\smallskip

Recall that the periodic solution of (\ref{eqB}) is $R$-symmetric if
and only if it intersects $\Fix(R)$ in exactly two points. In
conclusion, we obtain all small symmetric periodic solutions of
(\ref{eqB}) by solving the equation \begin{equation}\label{eqG}
G(u,\sigma)=B(u,\sigma)\mid_{\Fix(R)}=0. \end{equation}

\medskip

\section{Birkhoff normal form}

\smallskip

In this section we briefly discuss some points concerning the
Birkhoff normal form that will be useful in the sequel. The Birkhoff
normal form is useful because it preserves the symplectic structure.
In our cases if the vector field is in the Birkhoff normal form then
it is in the Belitskii normal form, and so we can apply Theorem
\ref{teowagen}.

\smallskip

The function $\{f,g\}=\omega(X_f,X_g)$ is called the {\it Poisson
bracket} of the smooth functions $f$ and $g$. Let ${\mathcal H}_n$
be the set of all homogeneous polynomials of degree $n$. The adjoint
map $Ad_{H_2} : {\mathcal H}_n \rightarrow {\mathcal H}_n$ is
defined by

\smallskip

\begin{equation}\label{adj} Ad_{H_2}(H) = \{H_2,H\} =
\omega(X_{H_2},X_H) = <-X_{H_2},\nabla H>.
\end{equation}

\smallskip

The {\em Birkhoff Normal Form Theorem} (cf. \cite{Ta,GH,Wi}) states
that if we have a Hamiltonian $H = H_2 + H_3 + H_4 + \cdots$,  where
$H_i \in {\mathcal H}_i$ is the homogeneous part of degree $i$, and
${\mathcal G}_i \subset {\mathcal H}_i$ satisfies ${\mathcal G}_i
\oplus Range(Ad_{H_2})={\mathcal H}_i$, then there exists a formal
symplectic power series transformation $\Phi$ such that $H\circ\Phi
= H_2 + \widetilde{H_3} + \widetilde{H_4} + \cdots$ where
$\widetilde{H_i} \in {\mathcal G}_i\ (i=3,4,\dots).$ In particular,
if $Ad_{H_2}$ is semi-simple, as in our case, then $Ker(Ad_{H_2})$
is the complement of $Range(Ad_{H_2})$.

\smallskip

As $R_j$ is symplectic, the change of coordinates $\Phi$ can be
chosen in such a way that $H\circ\Phi$ satisfies $H\circ\Phi\circ
R_j=-H\circ\Phi$. In order to see this, we can split ${\mathcal
H}_i={\mathcal H}_i^+ \oplus {\mathcal H}_i^-$, where ${\mathcal
H}_i^+=\{H\in{\mathcal H}_i:H\circ R_j=H\}$ and ${\mathcal
H}_i^-=\{H\in{\mathcal H}_i:H\circ R_j=-H\}$. If $R_j$ is
symplectic, then $Ad_{H_2}({\mathcal H}_i^\pm) = {\mathcal
H}_i^\mp$. In this case, if ${\mathcal H}_i = {\mathcal G}_i \oplus
Ad_{H_2}({\mathcal H}_i)$, then ${\mathcal H}_i^- = ({\mathcal
G}_i\cap{\mathcal H}_i^-) \oplus Ad_{H_2}({\mathcal H}_i^+)$. Now we
can perform the change of coordinates restricted to ${\mathcal
H}_i^-$. It implies that all monomial terms in the image of the
adjoint restricted to $H_i^-$ can be removed and it will remain only
monomials in the kernel of the adjoint restricted to $H_i^-$. And
so, the normal form is also $R_j$--reversible.

\begin{definition}\label{cond birk}
We say that a Hamiltonian vector field $X_H$ satisfies the Birkhoff
Condition (BC) if $Ad_{H_2}(H)=0$.
\end{definition}

\medskip

\begin{observation}
 By the equalities \eqref{adj}, the condition of the Definition
\ref{cond birk} is equivalent to $\omega(X_{H_2},X_H)=0$ or
$\{H_2,H\}=0$.
\end{observation}

\section{Two degrees of freedom}

\smallskip

In \cite{BT} a Birkhoff normal form for each $X\in \Omega^0$ is
derived and the following result is obtained:

\smallskip

\begin{theorem}
Assume $H$ is a Hamiltonian that is anti-invariant with respect to
the involution and the associated vector field $X_H$ has an
elliptical equilibrium point. Then there exists another Hamiltonian
$\widetilde{H}$, formally $C^k$--equivalent to $H$, such that the
vector field $X_{\widetilde{H}}$ has two one--parameter families of
symmetric periodic solutions, with period near $2\pi/\sqrt{ad-bc}$,
as in the Liapunov's Theorem, going through the equilibrium point.
\end{theorem}

\smallskip

Let $\Omega^0$ be the space of the $C^{\infty}$ $R_0$--reversible
Hamiltonian vector fields with two degrees of freedom in
$\mathbb{R}^4$ and fix the coordinate system $(x_1,y_1,x_2,y_2)$
$\in\mathbb{R}^4$. We endow $\Omega^0$ with the
$C^{\infty}$--topology. Let $\Omega^0_B\subset\Omega^0$ be the space
of the vector fields that satisfy the Birkhoff condition and
$\Omega^0_{\omega}\subset\Omega^0$ be the space of the analytic
ones. We prove the following result, which generalizes the previous
one.

\medskip

%\vspace{.5cm}

\noindent{\bf Theorem A:} {\it There exists an open set
$\mathcal{U}^0\subset\Omega^0_B$ (respec. $\Omega^0_{\omega}$) such
that
\begin{itemize}
\item[(a)]$\mathcal{U}^0$ is determined by the $3$--jet of the
vector fields.
\item[(b)] each $X\in\mathcal{U}^0$ possesses two $1$--parameter families of symmetric periodic
solutions terminating at the equilibrium point.
\end{itemize}}

\smallskip

\dem Fix on $\R^4$ a symplectic structure as in the Proposition
\ref{prop1c}. So the normal form of an involution has one of the
following form: $Id_{\mathbb{R}^4}$, or $-Id_{\mathbb{R}^4}$, or
$R_0=\left(\begin{array}{cccc}1&0&0&0\\0&1&0&0\\0&0&-1&0\\0&0&0&-1\end{array}\right)$.
We just work with $R_0$--reversible vector fields.

\smallskip

As in the cases in $\mathbb{R}^6$ we have that by the hypothesis the
Hamiltonian $H$ satisfies $H\circ R_0=-H$, so the linear part of the
vector field $\X$ is given by
\begin{equation}\label{lin_part}
A=\left(\begin{array}{cccc} 0&0&a&b\\0&0&c&d\\-d&b&0&0\\c&-a&0&0
  \end{array}\right),
\end{equation}
and their eigenvalues are $\{\pm\sqrt{bc-ad}, \pm\sqrt{bc-ad}\}.$ We
are interested in the case with $bc-ad<0.$ We call $\alpha =
\sqrt{ad-bc}$ and in order to obtain the Jordan canonical form of
the matrix $A$ we consider the transformation matrix
$$P=\left(\begin{array}{cccc}
0&\frac{-b}{\alpha}&0&\frac{-a}{\alpha}\\
0&\frac{-d}{\alpha}&0&\frac{-c}{\alpha}\\
0& 0& 1& 0\\1& 0& 0& 0
\end{array}\right).$$ After this transformation we obtain
$$\widehat{A}=P^{-1}.A.P=
\left(\begin{array}{cccccc}
0&\alpha&0&0\\
-\alpha&0&0&0\\
0&0&0&\alpha\\
0&0&-\alpha&0
\end{array}\right),$$ and $$\widehat{R_0}=P^{-1}.R_0.P=\left(\begin{array}{cccccc}
-1& 0& 0& 0\\
0& 1& 0& 0\\
0& 0& -1& 0\\
0& 0& 0& 1
\end{array}\right),$$ where $P^{-1}$ is the inverse matrix of $P$.

\smallskip

Performing a time rescaling we can assume that $\alpha =1$. We write
the canonical real Jordan form of $A$ as
$$\hat{A}= \displaystyle\left(\begin{array}{cccc}0&1&0&0\\-1&0&0&0\\
0&0&0&1\\0&0&-1&0\end{array}\right).$$

\smallskip

First we obtain the Belitskii normal form of $X_H$, by considering
$h:\mathbb{R}^4\rightarrow\mathbb{R}^4$ up to $3^{rd}$ order, which
is given by
$X_{H}(x_1,y_{1},x_{2},y_{2})=A[x_1,y_{1},x_{2},y_{2}]+h(x_1,y_{1},x_{2},y_{2})$;
and after we require the condition that the Belitskii normal form is
$\widehat{R_0}-$reversible, i. e, $X_{H}
\widehat{R_0}=-\widehat{R_0} X_{H}.$ Then the system obtained is
given by
\begin{equation}\label{belitskiiR4}\begin{array}{ll}
\dot{x_1}=&y_1+(e_{21}y_1+e_{23}y_2)(x_1^2+y_1^2)+e_{30}y_2(x_2^2+y_2^2)\\
&+(e_{16}x_1+e_{24}x_2)(y_1x_2-x_1y_2)+e_{26}y_2(x_1x_2+y_1y_2),\\
\dot{y_1}=&-x_1+(-e_{21}x_1-e_{23}x_2)(x_1^2+y_1^2)-e_{30}x_2(x_2^2+y_2^2)\\
&+(e_{16}y_1+e_{24}y_2)(y_1x_2-x_1y_2)-e_{26}x_2(x_1x_2+y_1y_2),\\
\dot{x_2}=&y_2+(-d_{15}y_1-d_{22}y_2)(x_1^2+y_1^2)-(d_{20}y_1+d_{29}y_2)(x_2^2+y_2^2)\\
&-(d_{17}y_1+d_{25}y_2)(x_1x_2+y_1y_2),\\
\dot{y_2}=&-x_2+(d_{15}x_1+d_{22}x_2)(x_1^2+y_1^2)+(d_{20}x_1+d_{29}x_2)(x_2^2+y_2^2)\\
&+(d_{17}x_1+d_{25}x_2)(x_1x_2+y_1y_2).
\end{array}\end{equation}

\smallskip

Now we use the fact that the vector field satisfies the Birkhoff
Condition. First of all we observe that the canonical symplectic
matrix
$$J=\left(\begin{array}{cccc}0 &1&
 0& 0\\-1& 0&  0& 0\\0& 0& 0& 1\\0& 0& -1& 0\end{array}\right),$$
after the linear change of coordinates $P$, is transformed into
\[
\widehat{J}=P^{-1}JP=\left(\begin{array}{cccc}0 &0&
 -1& 0\\0& 0&  0& -1\\1& 0& 0& 0\\0& 1& 0& 0\end{array}\right).
\]

\smallskip

We take a general Hamiltonian function $H:\R^4\rightarrow\R$ of
$4^{th}$ order, compute the kernel of $Ad_{H_2}$ defined on
\eqref{adj}, where $H_2$ is the homogeneous part of degree $2$ of
$H$, and require that $H$ satisfies $H\circ\widehat{R_0}=-H$. The
terms up to $3^{rd}$ order is given by
$h_b(x)=\widehat{J}\cdot\nabla H(x)$; its expression is

\smallskip

\begin{equation}\label{birkhoffR4}\begin{array}{ll}
\dot{x_1}=& y_1 + a_1y_1(x_1^2 + y_1^2) + a_2(2x_1x_2y_1 - x_1^2y_2
+y_1^2y_2) \\ & + a_3(3x_2^2y_1 - 2x_1x_2y_2 + y_1y_2^2),\\
\dot{y_1}=& -x_1 + (a_2y_1 + 2a_3y_2)(x_2y_1 - x_1y_2) -
x_1(a_1(x_1^2 + y_1^2)\\ & + a_2(x_1x_2 + y_1y_2) + a_3(x_2^2 + y_2^2)),\\
\dot{x_2}=& y_2 + (2a_1x_1 + a_2x_2)(-x_2y_1 + x_1y_2) +
y_2(a_1(x_1^2 + y_1^2) \\ &+ a_2(x_1x_2 + y_1y_2) + a_3(x_2^2 + y_2^2)),\\
\dot{y_2}=& -x_2 + (2a_1y_1 + a_2y_2)(-x_2y_1 + x_1y_2) -
x_2(a_1(x_1^2 + y_1^2) \\ &+ a_2(x_1x_2 + y_1y_2) + a_3(x_2^2 +
y_2^2)).
\end{array}\end{equation}

\smallskip

\begin{observation}
 We observe here that we can apply Theorem \ref{teowagen}, when
the vector field is in the Belitskii normal form. This is not a
restriction because if the vector field satisfies the Birkhoff
Condition then it is in the Belitskii Normal Form. It is easy to see
that if $\{H_2, H\}=0$ then $D(\{H_2, H\})=0$, and so $A_0^T X_H -
DX_H A_0^T(x)=0$. For example, in our case we have
\[\widehat{J}=\left(\begin{array}{cccc}0 &0&
 -1& 0\\0& 0&  0& -1\\1& 0& 0& 0\\0& 1& 0& 0\end{array}\right), \mbox{ } A_0=\left(\begin{array}{cccc}0&1&0&0\\-1&0&0&0\\
0&0&0&1\\0&0&-1&0\end{array}\right)\mbox{
 and }
 X_{H_2}=\left(\begin{array}{r}y_1\\-x_1\\y_2\\-x_2\end{array}\right).\]
The Birkhoff condition implies
$-y_1H_{x_1}+x_1H_{y_1}-y_2H_{x_2}+x_2H_{y_2} = 0$. So
\begin{equation}\label{birk}\begin{array}{r} H_{y_1}
-y_1H_{x_1x_1}+x_1H_{y_1x_1}-y_2H_{x_2x_1}+x_2H_{y_2x_1} =
0,\\
-H_{x_1} -y_1H_{x_1y_1}+x_1H_{y_1y_1}-y_2H_{x_2y_1}+x_2H_{y_2y_1} =
0,\\
H_{y_2} -y_1H_{x_1x_2}+x_1H_{y_1x_2}-y_2H_{x_2x_2}+x_2H_{y_2x_2} =
0,\\
-H_{x_2} -y_1H_{x_1y_2}+x_1H_{y_1y_2}-y_2H_{x_2y_2}+x_2H_{y_2y_2} =
0.
\end{array}\end{equation}
On the other hand if we compute $A_0^TX_H - DX_H A_0^T(x)$, we
obtain
\[ -\left(\begin{array}{c}-H_{y_2}\\H_{x_2}\\H_{y_1}\\-H_{x_1}\end{array}\right)
+ \left(\begin{array}{cccc}-H_{x_2x_1}&-H_{x_2y_1}&-H_{x_2x_2}&-H_{x_2y_2}\\
-H_{y_2x_1}&-H_{y_2y_1}&-H_{y_2x_2}&-H_{y_2y_2}\\
H_{x_1x_1}&H_{x_1y_1}&H_{x_1x_2}&H_{x_1y_2}\\
H_{y_1x_1}&H_{y_1y_1}&H_{y_1x_2}&H_{y_1y_2}\end{array}\right)\left(\begin{array}{c}y_1\\-x_1\\y_2\\-x_2\end{array}\right),\]
and by \eqref{birk} we have that $A_0^TX_H - DX_H A_0^T(x)=0$, i.e.
the system is in the Belitskii Normal Form.
\end{observation}

The Liapunov-Schmidt reduction gives us all small
$\widehat{R_0}$--symmetric periodic solutions by solving the
equation

\smallskip

\begin{equation*}\label{b2}B(x,\sigma)|_{x\in \Fix(\widehat{R_0})}=0,\end{equation*} with
\begin{equation}\label{b3}B(x,\sigma)=(1+\sigma)Sx-\hat{A} x-h_b(x),\,\,
x\in\mathbb{R}^4,\end{equation} where $S$ is the semi-simple part of
(unique) $S-N-$decomposition of $\hat{A}$. (See \cite{MM}).

\smallskip

In our case, $\hat{A}$ is semi-simple and
$\Fix(\widehat{R_0})=\{(0,y_{1},0,y_{2});\,\,y_1,y_{2}\in\mathbb{R}\}$.
Recall that the reduced equation, $B(x,\sigma)=0,$ is defined in
$\mathcal{N}\times \mathbb{R}$, where $\mathcal{N}=\{exp(\hat{A}
t)x; x\in V\}\in C^1_{2 \pi}$ and $V=span\{e_1,e_2,e_3,$ $e_4\}$.

The symplectic structure $\widehat{J}$ give us that $X_H$ is written
in the following form
$h_b(x)=h_b(x_1,y_1,x_2,y_2)=(-H_{x_2}(x_1,y_1,x_2,y_2),$
$-H_{y_2}(x_1,y_1,$ $x_2,y_2),$ $H_{x_1}(x_1,y_1,x_2,y_2),$
$H_{y_1}(x_1,y_1,x_2,y_2))$. Using the fact that $h_b$ satisfies the
Birkhoff Condition we have that
\[y_1H_{x_1}(x_1,y_1,x_2,y_2)-x_1H_{y_1}(x_1,y_1,x_2,y_2)+
y_2H_{x_2}(x_1,y_1,x_2,y_2)-\] \[x_2H_{y_2}(x_1,y_1,x_2,y_2)=0, \
\forall(x_1,y_1,x_2,y_2)\in\R^4.\] Hence at the points $(0,0,0,y_2)$
we have $y_2H_{x_2}(0,0,0,y_2)=0$. It implies that
$H_{x_2}(0,y_1,0,y_2)=y_1\bar{f}(y_1,y_2)$. Analogously we have that
$H_{x_1}(0,y_1,0,y_2)=y_2\bar{g}(y_1,y_2)$. So

\smallskip

\begin{equation}\label{equacao5}
G(y_{1},y_{2},\sigma)=B(x,\sigma)|_{x\in \Fix(\widehat{R_0})}=
\left[\begin{array}{l}
-y_1(a_1y_1^2 + a_2y_1y_2 + a_3y_2^2 - \sigma \mathbf{+\cdots})\\
-y_2(a_1y_1^2 + a_2y_1y_2 + a_3y_2^2 - \sigma\mathbf{+\cdots})
\end{array}\right].
\end{equation}

\smallskip

For the analytic case we have that the equation
$$G(y_{1},y_{2},\sigma)=(0,0)$$ is given by
\[\begin{array}{l}
-y_1(a_1y_1^2 + a_2y_1y_2 + a_3y_2^2 - \sigma) + H_1(y_1,y_2) = 0, \\
-y_2(a_1y_1^2 + a_2y_1y_2 + a_3y_2^2 - \sigma) + H_2(y_1,y_2) = 0,
\end{array}\]
and multiplying the first equation by $-y_2$ and the second by $y_1$
we get $y_2H_1=y_1H_2$. Using the fact that $H_1$ and $H_2$ are
analytic we have that there exists $\widetilde{H}$ such that
$H_1=y_1\widetilde{H}$ and $H_2=y_2\widetilde{H}$ for all
$(y_1,y_2)$.

If $a_1a_3\neq 0$ in \eqref{equacao5}, then we have two solutions
for the equation $G(y_1,y_2,\sigma)=0$. One solution is $y_1=0$ and
$y_2(\sigma)=\pm\sqrt{\frac{\sigma}{a_3}}+\dots$. And the second
solution is $y_2=0$ and
$y_1(\sigma)=\pm\sqrt{\frac{\sigma}{a_1}}+\dots$.

\smallskip

We define ${\mathcal U}^0=\mathcal{U}^0_1\cap\mathcal{U}^0_2$ where
$$\begin{array}{ll}
\mathcal{U}^0_1= & \left\{\begin{array}{ll}X\in\Omega^0_B;\,\,
&\mbox{the canonical form of $DX(0)$ satisfies $ad-bc>0$}\end{array}\right\}\\ \mbox{and} &\\
\mathcal{U}^0_2= & \left\{\begin{array}{ll}X\in\Omega^0_B;\,\,
&\mbox{the coefficients  of \eqref{birkhoffR4} satisfies $a_1a_3\neq
0$}
\end{array}\right\}.
\end{array}$$

\smallskip

In ${\mathcal U}^0=\mathcal{U}^0_1\cap\mathcal{U}^0_2
\subset\Omega^0_B$ for each $\sigma$ the equation
$G(y_{1},y_{2},\sigma)=0$ has two nonzero solutions terminating at
the origin when $\sigma$ is tending to $0$. So, in the original
problem we have two one parameter families of periodic solutions
terminating the origin (when $\sigma \rightarrow 0$). \cqd

\medskip

\section{Three degrees of freedom}

\smallskip

As in the previous section, let $\Omega^1$ (resp. $\Omega^2$) be the
space of the $C^{\infty}$ $R_1$-reversible (resp. $R_2$--reversible)
Hamiltonian vector fields with three degrees of freedom in
$\mathbb{R}^6$ and fix a coordinate system
$(x_1,y_1,x_2,y_2,x_3,y_3)\in\mathbb{R}^6$. We endow $\Omega^1$ and
$\Omega^2$ with the $C^{\infty}$--topology.  Let $\Omega^2_B$ (resp.
$\Omega^2_\omega$) be the space of vector fields in $\Omega^2$ that
satisfy the Birkhoff Condition (resp. that are analytic).

\medskip

\subsection{Case 6:2} $\;$

\medskip

\noindent {\bf Theorem B:} {\it There exists an open set
$\mathcal{U}^1\subset\Omega^1$ such that
\begin{itemize}
\item[(a)]$\mathcal{U}^1$ is determined by the $2$--jet of the
vector fields.
\item[(b)]for each $X\in\mathcal{U}^1$
there is no symmetric periodic orbit arbitrarily close to the
equilibrium point.
\end{itemize}}

\dem First we obtain the Belitskii normal form of $X_H$, by
considering $h:\mathbb{R}^6\rightarrow\mathbb{R}^6$ up to 2$^{nd}$
order, and then we require that the Belitskii normal form is
$\widehat{R_1}-$reversible, i. e, $X_{H}
\widehat{R_1}=-\widehat{R_1} X_{H}.$ After that we take the Birkhoff
normal form. The new symplectic structure is
$\widehat{J}=P_1^tJP_1$, where $P_1$ is the linear matrix that
brings the linear part of the vector field to the Jordan canonical
form. The Birkhoff normal form is obtained by taking a general
Hamiltonian function $H:\R^6\rightarrow\R$ of $3^{rd}$ order,
computing the kernel of $Ad_{H_2}$ and requiring that $H$ satisfies
$H\circ\widehat{R_1}=-H$. The Birkhoff normal form up to $2^{nd}$
order is given by $h_b(x)=\widehat{J}\cdot\nabla H(x)$. Finally, the
Liapunov-Schmidt reduction gives us all small
$\widehat{R_1}$--symmetric periodic solutions by solving the
equation
\begin{equation*}\label{b2}B(x,\sigma)|_{x\in \Fix(\widehat{R_1})}=0,\end{equation*} with
$$B(x,\sigma)=(1+\sigma)Sx-\widehat{A_1} x-h_b(x),\,\, x\in\mathbb{R}^6,$$
and $S$ is the semi-simple part of (unique) $S-N-$decomposition of
$\widehat{A_1}$. (See \cite{MM}). In our case, $\widehat{A_1}$ is
semi-simple and $\Fix(\widehat{R_1})=\{(0,0,x_{2},0,x_{3},0);\,\,$
$x_2,$ $x_{3}\in\mathbb{R}\}$. We recall that the reduced equation
of the Liapunov-Schmidt, $B(x,\sigma)=0,$ is defined in
$\mathcal{N}\times \mathbb{R}$, where
$\mathcal{N}=\{exp(\widehat{A_1} t)x; x\in V\}\in C^1_{2 \pi}$ and
$V=ger\{e_1,e_2,e_3,e_4,e_5,$ $e_6\}$.

We derive the following expression
\begin{equation}\label{b2}
\begin{array}{l}
G(x_{2},x_{3},\sigma)=B(x,\sigma)|_{x\in \Fix(\widehat{R_1})}= \\ \\
=\left[\begin{array}{l}
b_1x_2^2 + x_3(b_2x_2 + b_3x_3)+\cdots\\
b_4x_2^2 + x_3(b_5x_2 + b_6x_3)+\cdots\\
x_2(-\sigma  + \delta)+\cdots\\
x_3(-\sigma  + \delta)+\cdots
\end{array}\right].\end{array}\end{equation}

\smallskip

Observe that the equation $b_1x_2^2+b_2x_2x_3 + b_3x_3^2 =0$,
generically, either has the solution  $(x_2,x_3)=(0,0)$, or has a
pair of straight lines solutions given by
$(c_1x_2+d_1x_3)(c_2x_2+d_2x_3)=0$. The equation $b_4x_2^2+b_5x_2x_3
+ b_6x_3^2 =0$ is analogous. We can conclude that if the two first
components of (7.13) have no comom factor of the form $cx_2+dx_3$
then we have just the solution $(x_2,x_3)=(0,0)$ for the two
previous equations.

\smallskip

We define the following open sets:
$$\begin{array}{ll}
\mathcal{U}_1^1= & \left\{\begin{array}{ll}X\in\Omega^1;\,\,
&\mbox{the canonical form of $DX(0)$ satisfies \eqref{naosei} } \end{array}\right\},\\
\mathcal{U}_2^1=&\left\{\begin{array}{ll}X\in\Omega^1;\,\,
&\mbox{the $2$--jet of the two first equations of \eqref{b2} }\\ &
\mbox{ have no common factor }
\end{array}\right\}. \end{array}$$ Then
$\mathcal{U}^1=\mathcal{U}_1^1\cap \mathcal{U}_2^1$ is an open set
in $\Omega^1$.

\smallskip

The pair $(x_2,x_3)=(0,0)$ is the unique solution of the equation
$G=0$. So, near the origin there are no \textit{symmetric} periodic
orbits for this case. \cqd

\smallskip

\subsection{Case 6:4}$\;$

\medskip

\noindent {\bf Theorem C:} {\it There exists an open set
$\mathcal{U}^2\subset\Omega^2_B$ (respec. $\Omega^2_\omega$) such
that
\begin{itemize}
\item[(a)]$\mathcal{U}^2$ is determined by the $3$--jet of the
vector fields. \item[(b)]each $X\in\mathcal{U}^2$ has two
$2$--parameter families of periodic solutions
$\gamma_{\sigma,\lambda}^1$ and $\gamma_{\sigma,\lambda}^2$ with
$\sigma\in (-\epsilon,\epsilon)$ and $\lambda\in[0,2\pi]$, such
that, for each $\lambda_0$, $\lim_{\sigma\rightarrow
0}\gamma_{\sigma,\lambda_0}^j=0$, for $j=1,2$, and the periods tend to $2\pi/\alpha$ when $\sigma\rightarrow 0$. %In
%particular, if we choose $(r,s)=(\sigma,0)$ or $(r,s)=(0,\sigma)$
%in \eqref{familis} then we obtain two $1$--parameter Liapunov
%families for $X$.
\end{itemize}}

\smallskip

\dem First of all we derive the reversible Belitskii normal form of
$X_H$ up to $2^{nd}$ order. We observe that it coincides with the
reversible Birkhoff normal form and is given by:

\smallskip

\begin{equation}\label{birkhoff64}
X_{h_b}=\left[\begin{array}{l} -\dfrac{b (x_3 y_2-x_2 y_3)
\alpha ^2}{\beta }\\ \\
\dfrac{a (x_3 y_2-x_2y_3) \alpha ^2}{\beta }\\ \\
(-a x_1-b y_1) y_2+\alpha  y_2\\ \\
x_2 (a x_1+b y_1)-\alpha x_2  \\ \\
(-a x_1-b y_1) y_3+\alpha  y_3\\ \\
x_3 (a x_1+by_1)-\alpha x_3
\end{array}\right],\end{equation}
where $a=b_{65}/{\alpha}$, $b=b_{71}/{\alpha}$ e
${\alpha}=\sqrt{-a_{06} a_{10}+a_{05} a_{11}-a_{15} a_{17}+a_{14}
a_{18}}$.

\smallskip

As in the other cases, the Liapunov-Schmidt reduction gives us all
small $\widehat{R_2}$--symmetric periodic solutions by solving the
equation

\smallskip
\begin{equation*}\label{b3}B(x,\sigma)|_{x\in \Fix(\widehat{R_2})}=0,\end{equation*} with
$$B(x,\sigma)=(1+\sigma)Sx-\widehat{A_2} x-h_b(x),\,\, x\in\mathbb{R}^6.$$

\smallskip

As before $S$ is the semi-simple part of (the unique)
$S-N-$decomposition of $\widehat{A_2}$. (See \cite{MM}). In our
case, $\widehat{A_2}$ is semi-simple and
$\Fix(\widehat{R_2})=\{(x_1,y_1,0,y_{2},$ $0,$ $y_{3});$
$\,\,x_1,y_1,y_{2},y_{3}\in\mathbb{R}\}$. We recall that the reduced
equation of the Liapunov-Schmidt, $B(x,\sigma)=0,$ is defined on
$\mathcal{N}\times \mathbb{R}$, where
$\mathcal{N}=\{exp(\widehat{A_2} t)x;$ $x\in V\}\in C^1_{2 \pi}$ and
$V=ger\{e_1,e_2,e_3,e_4,e_5,$ $e_6\}$.

\smallskip

Like in the proof of Theorem A, we derive the following expression

\smallskip

\begin{equation}\label{b3}
\begin{array}{l}
G(x_1,y_1,y_{2},y_{3},\sigma)=B(x,\sigma)|_{x\in \Fix(\widehat{R_2})}= \\ \\
=\left[\begin{array}{l}
y_2(\sigma+a_{1}x_{1}+a_{2}y_{1} +a_3 x_1^2 + a_4y_1^2 + a_5 y_2^2 + a_6 y_2 y_3 + a_7 y_3^2+\cdots)\\
y_3(\sigma+a_{1}x_{1}+a_{2}y_{1}+a_3 x_1^2 + a_4y_1^2 + a_5 y_2^2 +
a_6 y_2 y_3 + a_7 y_3^2+\cdots)
\end{array}\right].\end{array}\end{equation}

\smallskip

If $a_5a_7\neq 0$ in \eqref{b3}, then for each $(x_1,y_1)$ close to
$(0,0)$ we have two solutions for the equation
$G(x_1,y_1,y_2,y_3,\sigma)=0$. One solution is $y_2=0$ and
$y_3(x_1,y_1,\sigma)=\pm\sqrt{\frac{\sigma+a_{1}x_{1}+a_{2}y_{1}}{a_7}}+\dots$.
And the second solution is $y_3=0$ and
$y_2(x_1,y_1,\sigma)=\pm\sqrt{\frac{\sigma+a_{1}x_{1}+a_{2}y_{1}}{a_5}}+\dots$.

We define the following open sets:
$$\begin{array}{ll}
\mathcal{U}_1^2= & \left\{\begin{array}{ll}X\in\Omega^2_B;\,\,
&\mbox{the canonical form of $DX(0)$ satisfies \eqref{naosei2} } \end{array}\right\},\\
\mathcal{U}_2^2=&\left\{\begin{array}{ll}X\in\Omega^2_B;\,\,
&\mbox{the coefficients of \eqref{b3} satisfies $a_5a_7\neq 0$ }
\end{array}\right\}. \end{array}$$ Then
$\mathcal{U}^2=\mathcal{U}_1^2\cap \mathcal{U}_2^2$ is an open set
in $\Omega^2_B$. For each $X\in\mathcal{U}^2$ and $\sigma$ we
consider $\gamma_{\sigma}^1:(x_1,y_1)\mapsto$ $(x_1,y_1,$ $0,$
$y_3(x_1,y_1,\sigma))$ and $\gamma_{\sigma}^2:(x_1,y_1)\mapsto$
$(x_1,y_1,$ $y_2(x_1,y_1,\sigma),$0$)$. Now we take the
parametrization $(x_1,y_1)\mapsto (a\sigma,b\sigma)$. We have
$\gamma_{\sigma\lambda_0}^1: (a\sigma,b\sigma)\mapsto
(a\sigma,b\sigma,0,$ $y_3(a\sigma,b\sigma,\sigma))$ and
$\gamma_{\sigma\lambda_0}^2: (a\sigma,b\sigma)\mapsto
(a\sigma,b\sigma,$ $y_2(a\sigma,b\sigma,\sigma),0)$ where
$\lambda_0= a/b$. Then, there exists two $2$--parameter family of
periodic orbits $\gamma_{\sigma\lambda}^1$ and
$\gamma_{\sigma\lambda}^2$ such that for each $\lambda_0\in\R$, the
families of periodic orbits $\gamma_{\sigma\lambda_0}^j$, for
$j=1,2$, are Liapunov families; i. e, $\lim_{\sigma\rightarrow
0}\gamma_{\sigma\lambda_0}^j=0$ and the period tends to
$2\pi/\alpha$. \cqd

\medskip

\section{Examples}

This section is devoted to present a mechanical example for the Case
$4:2$.

We consider two objects $m_1$ and $m_2$ with charge $q$ and $-q$.
They are at the position $(a,b)\in\R^2$ and $(-a,-b)\in\R^2$,
respectively. We assume that the system does not have kinetic
energy. So the total energy, i.e the Hamiltonian function is:
$$H(x,u,y,v)=\dfrac{-q}{\sqrt{(x-a)^2+(y-b)^2}}+\dfrac{q}{\sqrt{(x+a)^2+(y+b)^2}}.$$
Note that this Hamiltonian function satisfies the condition
$$H(\widehat{R}_0\cdot(x,u,y,v))=-H(x,u,y,v),$$ where
$\widehat{R}_0=\left(\begin{array}{cccccc}
-1& 0& 0& 0\\
0& 1& 0& 0\\
0& 0& -1& 0\\
0& 0& 0& 1
\end{array}\right).$

In another words, our system is a Hamiltonian
$\widehat{R}_0-$reversible vector field.

\begin{observation}
It is worth to say that the system \eqref{birkhoff64} (case 6 : 4 )
can be considered , in a similar way as \cite{YW}, a mathematical
model of a theoretical electrical circuit diagram.
\end{observation}

\noindent \textbf{Acknowledgements}\label{ackref} The authors thank
the dynamical system research group of Universitat Aut{\`o}noma de
Barcelona for the hospitality offered to us during part of the
preparation of this paper.
%\end{acknowledgements}

\end{document}